\documentstyle[amssymb,11pt]{amsart}
\makeatletter
\def\@cite#1#2{{\m@th\upshape\bfseries%
[{#1\if@tempswa{\m@th\upshape\mdseries, #2}\fi}]}}
\makeatother
\theoremstyle{plain}
\newtheorem{thm}{Theorem}[section]
\newtheorem{cor}[thm]{Corollary}
\newtheorem{prop}[thm]{Proposition}
\newtheorem{lem}[thm]{Lemma}
\theoremstyle{definition}
\newtheorem{rem}[thm]{Remark}
\newtheorem{defn}[thm]{Definition}

\renewcommand{\phi}{\varphi}
\newcommand{\upchi}{{\raise.35ex\hbox{$\chi$}}}
\def\<{\left<}
\def\>{\right>}
\def\cstar{$C^*$-algebra}
\begin{document}
\title{Generators of noncommutative dynamics}
%
\author{William Arveson}
\thanks{Partially supported by NSF grants DMS-9802474 and DMS-0100487} 
%
\address{Department of Mathematics,
University of California, Berkeley, CA 94720}
\email{arveson@@mail.math.berkeley.edu}

\begin{abstract}
For a fixed \cstar\ $A$, we consider all noncommutative 
dynamical systems that can be generated by $A$.  More 
precisely, an $A$-dynamical system is a triple $(i,B,\alpha)$ 
where $\alpha$ is a $*$-endomorphism of a \cstar\ $B$, 
and $i: A\subseteq B$ is the inclusion of $A$ as a 
$C^*$-subalgebra with the property that $B$ is 
generated by $A\cup \alpha(A)\cup \alpha^2(A)\cup\cdots$.  
There is a natural hierarchy in the class of $A$-dynamical 
systems, and there is a universal one that dominates 
all  others, denoted $(i,{\mathcal P}A,\alpha)$.  
We establish certain properties of $(i,{\mathcal P}A,\alpha)$ 
and give applications to some concrete issues of 
noncommutative dynamics.

For example, we show that every contractive completely 
positive linear map $\phi: A\to A$ gives rise to 
to a unique $A$-dynamical system $(i,B,\alpha)$ that 
is ``minimal" with respect to $\phi$, and 
we show that its \cstar\ $B$ 
can be embedded in the multiplier 
algebra of $A\otimes {\mathcal K}$.  
\end{abstract}
\maketitle

\section{Generators}

The flow of time in quantum theory is represented by
a one-parameter group of $*$-automorphisms $\{\alpha_t: t\in\mathbb R\}$
of a \cstar\ $B$.  There is often a $C^*$-subalgebra $A\subseteq B$ 
that can be singled out from physical considerations 
which, together with its time 
translates, generates $B$.  For example, in nonrelativistic
quantum mechanics the flow of time is represented by a 
one-parameter group of automorphisms of ${\mathcal B}(H)$, and the 
set of all 
bounded continuous functions of the configuration 
observables at time $0$ is a commutative \cstar\ $A$.   
The set of all time translates $\alpha_t(A)$ of $A$ 
generates an irreducible $C^*$-subalgebra 
$B$ of ${\mathcal B}(H)$.  In particular, for 
different times $t_1\neq t_2$, the \cstar s $\alpha_{t_1}(A)$ 
and $\alpha_{t_2}(A)$ do not commute with each other. Indeed, 
no nontrivial relations appear to exist between 
$\alpha_{t_1}(A)$ and $\alpha_{t_2}(A)$
when $t_1\neq t_2$.  

In this paper we look closely at this phenomenon, in 
a simpler but analogous setting.  Let $A$ be a \cstar, fixed 
throught. 
\begin{defn}\label{dynSystemDef}
An $A$-dynamical system is a triple $(i,B,\alpha)$ 
consisting of a $*$-endomorphism $\alpha$ acting 
on a \cstar\ $B$ and an injective $*$-homomorphism 
$i: A\to B$, such that $B$ is generated 
by $i(A)\cup\alpha(i(A))\cup\alpha^2(i(A))\cup\cdots$.  
\end{defn}
\noindent
We lighten notation by identifying $A$ with its image 
$i(A)$ in $B$, thereby replacing $i$ with the inclusion 
map $i: A\subseteq B$.  Thus, an $A$-dynamical system 
is a dynamical system $(B,\alpha)$ that contains $A$ as a 
$C^*$-subalgebra in a specified way, with the property that 
$B$ is the norm-closed linear span of finite 
products of the following form
\begin{equation}\label{spanEq}
B=\overline{\rm span}\{\alpha^{n_1}(a_1)\alpha^{n_2}(a_2)\cdots\alpha^{n_k}(a_k)\}
\end{equation}
where $n_1,\dots,n_k\geq 0$, $a_1,\dots,a_k\in A$, $k=1,2,\dots$

Our aim is to say something sensible about the class 
of {\it all} $A$-dynamical systems, and to obtain more 
detailed information about certain of its members.  
The opening paragraph illustrates the fact that in 
even the simplest cases, 
where $A$ is $C(X)$ or even a matrix algebra, 
the structure of individual 
$A$-dynamical systems can be very complex.   

There is a natural hierarchy in the class of all 
$A$-dynamical systems, defined by 
$(i_1,B_1,\alpha_1)\geq (i_2,B_2,\alpha_2)$ iff there 
is a $*$-homomorphism $\theta: B_1\to B_2$ satisfying 
$\theta\circ\alpha_1=\alpha_2\circ\theta$ and 
$\theta(a)=a$ for $a\in A$.  Since $\theta$ fixes $A$, 
it follows from (\ref{spanEq}) that  $\theta$ must be 
surjective, $\theta(B_1)=B_2$, hence 
$(i_2,B_2,\alpha_2)$ is a {\it quotient} of $(i_1,B_1,\alpha_1)$.  
Two $A$-dynamical systems are said to be {\it equivalent} if there is 
a map $\theta$ as above that is an isomorphism of \cstar s.  
This will be the case iff each of the $A$-dynamical systems 
dominates the other.  
One may also think of the class 
of all $A$-dynamical systems as a category, whose objects 
are $A$-dynamical systems and whose maps 
$\theta$ are described above.  

There is a largest equivalence class in this hierarchy, 
whose representatives are called 
{\it universal} $A$-dynamical systems.  We exhibit one
as follows.  
Consider the free product 
of an infinite sequence of copies of $A$, 
$$
{\mathcal P} A=A*A*\cdots.  
$$
Thus, we have a sequence of $*$-homomorhisms $\theta_0, \theta_1,\dots$ 
of $A$ into the \cstar\ ${\mathcal P}A$ such ${\mathcal P}A$ is generated 
by $\theta_0(A)\cup\theta_1(A)\cup\cdots$ and such that the following 
universal property is satisfied: for every  sequence 
$\pi_0,\pi_1,\dots$ of $*$-homomorphisms of $A$ into 
some other \cstar\ $B$, there is a unique  
$*$-homomorphism $\rho: {\mathcal P}A\to B$ such that 
$\pi_k=\rho\circ\theta_k$, $k=0,1,\dots$.  
Nondegenerate representations 
of ${\mathcal P}A$ correspond to sequences $\bar\pi=(\pi_0,\pi_1,\dots)$ 
of representations $\pi_k: A\to {\mathcal B}(H)$ of $A$ on 
a common Hilbert space $H$, subject to no condition other 
than the triviality of their common nullspace 
$$
\xi\in H,\quad \pi_k(A)\xi=\{0\}, \quad k=0,1,\dots \implies \xi=0.  
$$
A simple argument establishes the existence of ${\mathcal P}A$ by taking 
the direct sum of a sufficiently large set of 
such representation sequences $\bar\pi$.  

This definition does not exhibit ${\mathcal P}A$ in concrete 
terms (see \S \ref{S:construction} for that), 
but it does allow us to define a  
universal $A$-dynamical system.  
The universal property of ${\mathcal P}A$ implies that there is 
a shift endomorphism $\sigma: {\mathcal P}A\to {\mathcal P}A$ 
defined uniquely by $\sigma\circ\theta_k=\theta_{k+1}$, 
$k=0,1,\dots$.  It is quite easy to verify 
that $\theta_0$ is an injective $*$-homomorphism 
of $A$ in ${\mathcal P}A$,  and we use this map to 
identify $A$ with $\theta_0(A)\subseteq{\mathcal P}A$.  
Thus {\it the triple $(i,{\mathcal P}A,\sigma)$ 
becomes an $A$-dynamical system with the property that every 
other $A$-dynamical system is subordinate to it}.

Before introducing $\alpha$-expectations, we 
review some common terminology \cite{gPedAutbook}.  
Let $A\subseteq B$ be an inclusion of \cstar s. 
For any subset $S$ of $B$ we write $[S]$ for 
the norm-closed linear span of $S$.  
The subalgebra 
$A$ is said to be {\it essential} if the 
two-sided ideal $[BAB]$ it generates is 
an essential ideal 
$$
x\in B, \quad xBAB=\{0\}\implies x=0.  
$$
It is called {\it hereditary} 
if for $a\in A$ and $b\in B$, one has 
$$
0\leq b\leq a \implies b\in A.  
$$
The hereditary subalgebra of $B$ generated by 
a subalgebra $A$ is the 
closed linear span $[ABA]$ of all products $axb$, 
$a,b\in A$, $x\in B$, 
and in general 
$A\subseteq [ABA]$.  A {\it corner} of $B$ is a hereditary 
subalgebra of the particular 
form $A=pBp$ where $p$ is a projection 
in the multiplier algebra $M(B)$ of $B$.  

We also make essential 
use of {\it conditional expectations} 
$E: B\to A$.  A conditional expectation 
is an idempotent positive linear map 
with  range $A$, satisfying $E(ax)=aE(x)$
for $a\in A$, $x\in B$.  
When $A=pBp$ is a {\it corner} 
of $B$, the map 
$E(x)=pxp$, defines a conditional 
expectation of $B$ onto $A$.  On the other hand, many of
the  conditional expectations encountered 
here do not have this simple 
form, even when $A$ has a unit.  Indeed, 
if $A$ is subalgebra of $B$ that is 
{\it not} hereditary, then there is no natural conditional 
expectation $E: B\to A$.  In general, conditional expectations 
are completely positive linear maps with $\|E\|=1$.  

\begin{defn}\label{alphaDef} Let 
$(i,B,\alpha)$ be an $A$-dynamical system.  An 
$\alpha$-expectation is a conditional 
expectation $E: B\to A$ having the following 
two properties:
\begin{enumerate}
\item[E1.]
Equivariance:  $E\circ\alpha=E\circ\alpha\circ E$.  
\item[E2.]
The restriction of $E$ to the hereditary subalgebra 
generated by $A$ is multiplicative, 
$E(xy)=E(x)E(y)$, $x,y\in [ABA]$.  
\end{enumerate}
\end{defn}

Note that an {\it arbitrary} 
conditional expectation $E: B\to A$ 
gives rise to a linear map $\phi: A\to A$ by way 
of  $\phi(a)=E(\alpha(a))$, 
$a\in A$.  Such a 
$\phi$ is a completely positive map 
satisfying $\|\phi\|\leq 1$.   
Axiom E1 makes the assertion 
\begin{equation}\label{equivariance}
E\circ\alpha=\phi\circ E.   
\end{equation}
where $\phi=E\circ\alpha\upharpoonright_A$ is 
the linear map of $A$ associated with $E$.  

Property E2 is of course automatic 
if $A$ is a hereditary subalgebra of $B$.  
It is a fundamentally {\it noncommutative} 
hypothesis on $B$.  For example, if $Y$ is 
a compact Hausdorff space and $B=C(Y)$, then every 
unital subalgebra $A\subseteq C(Y)$ generates 
$C(Y)$ as a hereditary algebra.  
Thus the only linear maps 
$E: C(Y)\to A$ satisfying E2 
are $*$-endomorphisms of $C(Y)$.  
The key property of the universal $A$-dynamical system 
$(i,{\mathcal P}A,\sigma)$ follows.

\begin{thm}\label{universalDilationThm}
For every completely positive 
contraction $\phi: A\to A$, there 
is a unique $\sigma$-expectation $E: {\mathcal P}A\to A$ 
satisfying
\begin{equation}\label{compressionEq}
\phi(a)=E(\sigma(a)), \qquad a\in A.  
\end{equation}
\end{thm}

Both assertions
are nontrivial.  We prove 
uniqueness in the following section, 
see Theorem \ref{uniquenessThm}.  
Existence is taken up in \S \ref{S:construction}, 
see Theorem \ref{cpThm}.  

\section{moment polynomials}\label{S:moments}
This theory of generators rests on properties 
of certain noncommutative polynomials that are  
defined recursively as follows.  

\begin{prop}\label{momPolyProp}
Let $A$ be an algebra over a field $\mathbb F$.  For every 
linear map $\phi: A\to A$, 
there is a unique sequence of multilinear mappings 
from $A$ to itself, indexed by the 
$k$-tuples of nonnegative integers, 
$k= 1,2,\dots$, where for a 
fixed $k$-tuple $\bar n=(n_1,\dots,n_k)$
$$
a_1,\dots,a_k\in A\mapsto [\bar n; a_1, \dots,a_k]\in A
$$
is a $k$-linear mapping, all 
of which satisfy 
\begin{enumerate}
\item[MP1.]\label{mp1}
$\phi([\bar n; a_1,\dots,a_k])=[n_1+1,n_2+1,\dots,n_k+1;
a_1,\dots,a_k]$.  
\item[MP2.]\label{mp2}
Given a $k$-tuple for which $n_\ell=0$ for some 
$\ell$ between $1$ and $k$, 
$$
[\bar n; a_1,\dots,a_k]=[n_1,\dots,n_{\ell-1};
a_1,\dots,a_{\ell-1}]a_\ell 
[n_{\ell+1},\dots,n_k; a_{\ell+1},\dots,a_k].
$$
\end{enumerate}
\end{prop}

\begin{rem}
The proofs of both existence and uniqueness are straightforward 
arguments using induction on the number $k$ of variables, and 
we omit them.  
Note that in the second axiom MP2, 
we make the natural conventions when 
$\ell$  has one of the extreme values $1$, $k$.  For example, if  
$\ell=1$, then MP2 should be interpreted as 
$$
[0,n_2,\dots,n_k; a_1,\dots,a_k]=a_1[n_2,\dots,n_k; a_2,\dots,a_k].  
$$
In particular, in the linear case $k=1$, MP2
makes the assertion
$$
[0;a]=a, \qquad a\in A;
$$
and after repeated applications of axiom MP1 one obtains 
$$
[n;a]=\phi^n(a), \qquad a\in A, \quad n=0,1,\dots.  
$$
One may calculate any particular
moment polynomial explicitly, but 
the computations quickly become a tedious exercise 
in the arrangement of parentheses.  For example, 
\begin{align*}
[2,6,3,4;a,b,c,d]&=\phi^2(a\phi(\phi^3(b)c\phi(d))), \\
 [6,4,2,3;a,b,c,d]&=\phi^2(\phi^2(\phi^2(a)b)c\phi(d)).  
\end{align*}

Finally, we remark that when $A$ is a $C^*$-algebra 
and $\phi: A\to A$ is a 
linear map satisfying $\phi(a)^*=\phi(a^*)$, 
$a\in A$, then its associated moment polynomials 
obey the following symmetry
\begin{equation}\label{symmetryEq}
[n_1,\dots,n_k;a_1,\dots,a_k]^*=[n_k,\dots,n_1;a_k^*,\dots,a_1^*].  
\end{equation}
Indeed, one finds that the 
sequence of polynomials $[[\cdot;\cdot]]$ 
defined by 
$$
[[n_1,\dots,n_k;a_1,\dots,a_k]]=[n_k,\dots,n_1;a_k^*,\dots,a_1^*]^*
$$
also satisfies axioms MP1 and MP2, and hence must coincide 
with the moment polynomials of $\phi$ by the uniqueness 
assertion of Proposition \ref{momPolyProp}.  
\end{rem}

These polynomials are important because they are the 
expectation values of certain $A$-dynamical systems.  
\begin{thm}\label{uniquenessThm}
Let $\phi: A\to A$ be a completely 
positive map on $A$, satisfying $\|\phi\|\leq 1$, with associated 
moment polynomials $[n_1,\dots,n_k;a_1,\dots,a_k]$.  

Let $(i,B,\alpha)$ be an $A$-dynamical system 
and let $E: B\to A$ be an $\alpha$-expectation with the property 
$E(\alpha(a))=\phi(a)$, $a\in A$.  Then 
\begin{equation}
E(\alpha^{n_1}(a_1)\alpha^{n_2}(a_2)\cdots\alpha^{n_k}(a_k)) 
=[n_1,\dots,n_k; a_1,\dots,a_k].  \label{nPointEq}
\end{equation}
for every $k=1,2,\dots$, $n_k\geq 0$, $a_k\in A$.  
In particular, there is at most one 
$\alpha$-expectation $E: B\to A$ 
satisfying $E(\alpha(a))=\phi(a)$, $a\in A$.  
\end{thm}

\begin{Prf} One applies 
the uniqueness of moment polynomials as follows.  
Properties 
E1 and E2 of Definition \ref{alphaDef} imply that the 
sequence of polynomials $[[\cdot;\cdot]]$ defined by 
$$
[[n_1,\dots,n_k;a_1,\dots,a_k]]=E(\alpha^{n_1}(a_1)\cdots\alpha^{n_k}(a_k))
$$
must satisfy the two axioms MP1 and MP2.  Notice here that
E2 implies 
\begin{equation}\label{multipEq}
E(xay)=E(x)aE(y)\qquad x,y\in B, \quad a\in A,
\end{equation} 
since for an approximate 
unit $e_n$ for $A$ we can write $E(xay)$ as follows:
\begin{align*}
\lim_{n\to\infty}e_nE(xae_ny)e_n&=
\lim_{n\to \infty}E(e_nxae_nye_n)
=\lim_{n\to\infty}E(e_nxa)E(e_nye_n)\\
&=\lim_{n\to\infty}e_nE(x)ae_nE(y)e_n=E(x)aE(y).  
\end{align*}
Thus formula 
(\ref{nPointEq}) follows from the uniqueness assertion 
of Proposition \ref{momPolyProp}.  
The uniqueness of the $\alpha$-expectation 
associated with $\phi$ 
is now apparent from formulas (\ref{nPointEq}) and 
(\ref{spanEq}).  \qed \end{Prf}

\section{Existence of $\alpha$-expectations}\label{S:construction}
In this section we show that every completely positive map 
$\phi: A\to A$, with $\|\phi\|\leq 1$, 
gives rise to an $\alpha$-expectation 
$E: {\mathcal P}A\to A$ that is related to 
the moment polynomials of $\phi$ as in (\ref{nPointEq}).    
This is established through a construction 
that exhibits ${\mathcal P}A$ as 
the enveloping \cstar\ of a Banach $*$-algebra $\ell^1(\Sigma)$, 
in such a way that the desired conditional 
expectation appears as a completely positive map on 
$\ell^1(\Sigma)$.  The details are as follows.  

Let $S$ be 
the set of finite sequences $\bar n=(n_1,n_2,\dots,n_k)$
of nonnegative integers, $k=1,2,\dots$ 
which have distinct neighbors, 
$$
n_1\neq n_2, n_2\neq n_3,\dots, n_{k-1}\neq n_k.  
$$
Multiplication and  
involution are defined in $S$ as follows.  
The product of two elements 
$\bar m=(m_1,\dots,m_k), \bar n = (n_1,\dots,n_\ell)\in S$
is defined by conditional concatenation 
$$
\bar m\cdot \bar n=
\cases
(m_1,\dots,m_k, n_1,\dots, n_\ell), \qquad\text{if } m_k\neq n_1,\\
(m_1,\dots,m_k, n_2,\dots, n_\ell), \qquad\text{if } m_k= n_1,
\endcases
$$
where we make the natural conventions when $\bar n = (q)$ 
is of length 1, namely $\bar m\cdot (q)=(m_1,\dots, m_k,q)$ 
if $m_k\neq q$, and $\bar m\cdot (q)=\bar m$ 
if $m_k=q$.  
The involution in $S$ is defined by reversing the 
order of components  
$$
(m_1,\dots, m_k)^*=(m_k,\dots,m_1).  
$$
One  finds that $S$ is an associative $*$-semigroup.

Fixing a \cstar\ $A$, we 
attach a Banach space $\Sigma_\nu$ to  
every $k$-tuple $\nu=(n_1,\dots,n_k)\in S$ as follows
$$
\Sigma_\nu = 
\underbrace{A\hat\otimes\cdots\hat\otimes A}_{k\text{\ times}},
$$
the $k$-fold projective tensor product of copies of the 
Banach space $A$.  We assemble the $\Sigma_\nu$ 
into a family of Banach spaces 
over $S$, $p: \Sigma\to S$, by way of 
$\Sigma = \{(\nu, \xi): \nu\in S, \xi\in E_\nu\}$, $p(\nu,\xi)=\nu$.  

We introduce a multiplication 
in $\Sigma$ as follows.  Fix 
$\mu=(m_1,\dots, m_k)$ and $\nu=(n_1,\dots,n_\ell)$ in $S$ 
and choose 
$\xi\in \Sigma_\mu$, $\eta\in \Sigma_\nu$.  If $m_k\neq n_1$ 
then $\xi\cdot\eta$ is defined as the tensor product 
$\xi\otimes\eta \in \Sigma_{\mu\cdot\nu}$.  If 
$m_k=n_1$ then we must tensor over $A$ and make the 
obvious identifications.  More explicitly, in this case 
there is a natural map of 
the tensor product $\Sigma_\mu\otimes_A \Sigma_\nu$ 
onto 
$\Sigma_{\mu\cdot\nu}$ by making identifications 
of elementary tensors as follows:
$$
(a_1\otimes\cdots\otimes a_k)\otimes_A(b_1\otimes\cdots\otimes b_\ell)
\sim
a_1\otimes\cdots\otimes a_{k-1}\otimes 
a_kb_1\otimes b_2\otimes\cdots\otimes b_\ell.  
$$
With this convention  
$\xi\cdot\eta$ is defined by
$$
\xi\cdot\eta=\xi\otimes_A\eta\in \Sigma_{\mu\cdot\nu}.  
$$ 
This defines an associative multiplication in the family 
of Banach spaces $\Sigma$.  There is also a natural 
involution in $\Sigma$, defined on each 
$\Sigma_\mu$, $\mu=(m_1,\dots,m_k)$ as the unique 
antilinear isometry to $\Sigma_{\mu^*}$ satisfying
$$
((m_1,\dots,m_k), a_1\otimes\cdots\otimes a_k)^* =
((m_k,\dots,m_1), a_k^*\otimes\cdots\otimes a_1^*).  
$$
This defines an isometric antilinear mapping of the 
Banach space $\Sigma_\mu$ onto $\Sigma_{\mu^*}$, for 
each $\mu\in S$, and thus the structure $\Sigma$ becomes 
an involutive $*$-semigroup in which each fiber 
$\Sigma_\mu$ is a Banach space.  

Let $\ell^1(\Sigma)$ be the Banach $*$-algebra 
of summable sections.  
The norm and involution are the natural ones
$\|f\|=\sum_{\mu\in \Sigma}\|f(\mu)\|$, 
$f^*(\mu)=f(\mu^*)^*$.  Noting that 
$\Sigma_\lambda\cdot\Sigma_\mu\subseteq\Sigma_{\lambda\cdot\mu}$, 
the multiplication in
$\ell^1(\Sigma)$ is defined by convolution
$$
f*g(\nu)=\sum_{\lambda\cdot\mu=\nu}f(\lambda)\cdot g(\mu), 
$$
and one easily verifies that $\ell^1(\Sigma)$ is a Banach 
$*$-algebra. 

For $\mu=(m_1,\dots,m_k)\in S$ and $a_1,\dots,a_k\in A$ we 
define the function 
$$
\delta_\mu\cdot a_1\otimes\dots\otimes a_k\in \ell^1(\Sigma)
$$
to be zero except at $\mu$, and at $\mu$ it 
has the value $a_1\otimes\dots\otimes a_k\in \Sigma_\mu$.  
These elementary functions have $\ell^1(\Sigma)$ as their 
closed linear span.  Finally, there is a natural sequence of 
$*$-homomorphisms 
$\theta_0, \theta_1,\dots: A\to \ell^1(\Sigma)$ defined by 
$$
\theta_k(a)=\delta_{(k)}\cdot a, \qquad a\in A,\quad k=0,1,\dots, 
$$
and these maps are related to the generating sections by
$$
\delta_{(n_1,\dots,n_k)}\cdot a_1\otimes\dots\otimes a_k=
\theta_{n_1}(a_1)\theta_{n_2}(a_2)\cdots\theta_{n_k}(a_k).  
$$

The algebra $\ell^1(\Sigma)$ fails to have a unit, but it has the 
same representation theory 
as ${\mathcal P}A$ in the following sense.  
Given a sequence of representations 
$\pi_k:A\to{\mathcal B}(H)$, $k=0,1,\dots$, 
fix $\nu=(n_1,\dots,n_k)\in S$.
There is a unique bounded linear operator   
$L_\nu: \Sigma_\nu\to{\mathcal B}(H)$ of norm 
$1$ that is defined by its action on elementary 
tensors as follows 
$$
L_\nu( a_1\otimes\cdots\otimes a_k)=
\pi_{n_1}(a_1)\cdots\pi_{n_k}(a_k).  
$$
Thus there is a bounded 
linear map $\tilde\pi: \ell^1(\Sigma)\to{\mathcal B}(H)$ 
defined by 
$$
\tilde\pi(f)=\sum_{\mu\in S}L_\mu(f(\mu)), 
\qquad f\in \ell^1(\Sigma).     
$$
One finds that $\tilde\pi$ 
is a $*$-representation of 
$\ell^1(\Sigma)$ with
$\|\tilde\pi\| = 1$.  
This representation satisfies 
$\tilde\pi\circ\theta_k=\pi_k$, 
$k=0,1,2,,\dots$.  Conversely, every bounded 
$*$-representation $\tilde\pi$ of $\ell^1(\Sigma)$ 
on a Hilbert space $H$ is associated with 
a sequence of representations 
$\pi_0, \pi_1,\dots$ of $A$ on $H$ by way
of $\pi_k=\tilde\pi\circ\theta_k$.

The results of the preceding discussion 
are summarized as follows: 

\begin{prop}
The enveloping \cstar\ $C^*(\ell^1(\Sigma))$, together 
with the sequence of homomorphisms 
$\tilde\theta_0, \tilde\theta_1,\dots: A\to C^*(\ell^1(\Sigma))$ 
defined by the maps 
$\theta_0, \theta_1,\dots: A\to \ell^1(\Sigma)$, has the 
same universal property as the infinite free product 
${\mathcal P}A=A*A*\cdots$, and is therefore isomorphic 
to ${\mathcal P}A$.  
\end{prop}

Notice that the natural shift endomorphism of $\ell^1(\Sigma)$ 
is defined by 
$$
\sigma:\delta_{(n_1,\dots,n_k)}\cdot\xi\mapsto
\delta_{(n_1+1,\dots,n_k+1)}\cdot\xi, \qquad 
\nu=(n_1,\dots,n_k)\in \Sigma,\quad\xi\in \Sigma_\nu
$$ 
and it promotes to the natural 
shift endomorphism of ${\mathcal P}A=C^*(\ell^1(\Sigma))$.  
The inclusion of $A$ in $\ell^1(\Sigma)$ is given by 
the map $\theta_0(a)= \delta_{(0)}a\in \ell^1(\Sigma)$, 
and it too promotes 
to the natural inclusion of $A$ in ${\mathcal P}A$.  

Finally, we fix a contractive completely 
positive map $\phi: A\to A$, and consider the moment 
polynomials associated with it by Proposition \ref{momPolyProp}.  
A straightforward argument shows that there is a unique 
bounded linear map $E_0: \ell^1(\Sigma)\to A$ 
satisfying 
$$
E_0(\delta_{(n_1,\dots,n_k)}\cdot a_1\otimes\cdots\otimes a_k) 
=[n_1,\dots,n_k; a_1,\dots,a_k], 
$$
for $(n_1,\dots,n_k)\in S$, $a_1,\dots,a_k\in A$, 
$k=1,2,\dots$, and $\|E_0\|= \|\phi\|\leq 1$.  
Using the axioms MP1 and 
MP2, one finds that the map 
$E_0$ preserves the adjoint (see Equation (\ref{symmetryEq})),  
satisfies the conditional expectation property 
$E_0(af)=aE_0(f)$ for $a\in A$, $f\in\ell^1(\Sigma)$, 
that the restriction of $E_0$ to the 
``hereditary" $*$-subalgebra 
of $\ell^1(\Sigma)$ spanned by 
$\theta_0(A)\ell^1(\Sigma)\theta_0(A)$ is multiplicative, 
and that it is related to $\phi$ 
by $E_0\circ\sigma=\phi\circ E_0$ and 
$E_0(\sigma(a))=\phi(a)$, $a\in A$.  Thus, $E_0$ 
satisfies the axioms of Definition \ref{alphaDef}, suitably 
interpreted for the Banach $*$-algebra $\ell^1(\Sigma)$.  

In view of the basic fact that a bounded completely positive 
linear map of a Banach $*$-algebra to $A$ 
promotes naturally to a completely positive map of 
its enveloping \cstar \ to $A$, the critical 
property of $E_0$ reduces to:

\begin{thm}\label{cpThm}
For every $n\geq 1$,  
$a_1,\dots,a_n\in A$, and $f_1, \dots,f_n\in\ell^1(\Sigma)$, 
we have 
$$
\sum_{i,j=1}^n a_j^*E_0(f_j^*f_i)a_i\geq 0.  
$$
Consequently, $E_0$ extends uniquely through the 
completion map $\ell^1(\Sigma)\to{\mathcal P}A$ to 
a completely positive map $E_\phi: {\mathcal P}A\to A$ that 
becomes a $\sigma$-expectation 
satisfying Equation (\ref{compressionEq}).  
\end{thm}

We sketch the proof of Theorem \ref{cpThm}, detailing the 
critical steps.   Using the fact that $\ell^1(\Sigma)$ 
is spanned by the generating family 
$$
G=\{\delta_{(n_1,\dots,n_k)}\cdot a_1\otimes\dots\otimes a_k:
(n_1,\dots,n_k)\in S, \quad a_1,\dots, a_k\in A, \quad k\geq 1\}
$$
one easily reduces the proof of Theorem \ref{cpThm} 
to the following more concrete assertion: for 
any finite set of elements $u_1,\dots,u_n$ in 
$G$, 
the $n\times n$ matrix  $(a_{ij})=(E_0(u_j^*u_i))\in M_n(A)$ 
is positive.  

The latter is established by an inductive 
argument on the ``maximum height" 
$\max( h(u_1),\dots,h(u_n))$, where 
the height of an element 
$u=\delta_{(n_1,\dots,n_k)}\cdot a_1\otimes\cdots\otimes a_k$
in $G$ is defined as $h(u)=\max(n_1,\dots,n_k)$.   
The general case easily reduces to that 
in which $A$ has a unit $e$, and in that setting 
the inductive step is implemented by the following.  

\begin{lem}\label{keyLemma}
Choose $u_1,\dots,u_n\in G$ such that the maximum height 
$N=\max(h(u_1),\dots,h(u_n))$ is positive.  For 
$k=1,\dots,n$ there are elements $b_k, c_k\in A$ and 
$v_k\in G$ such that $h(v_k)<N$ and 
\begin{equation}\label{keyEq}
E_0(u_j^*u_i)=b_j^*\phi(E_0(v_j^*v_i))b_i
+c_j^*(e-\phi(e))c_i, \qquad 1\leq i,j\leq n.  
\end{equation}
\end{lem}

\begin{rem}
Note that if an inductive hypothesis provides positive 
$n\times n$ matrices of the form 
$(E_0(v_j^*v_i))$ whenever $v_1,\dots,v_n\in G$ 
have height $<N$,  
then the $n\times n$ matrix whose $ij$th term is the 
right side of (\ref{keyEq}) must also be positive,  
because $\phi$ is a completely positive 
map and $0\leq \phi(e)\leq e$.
It follows from Lemma \ref{keyLemma} that 
$(E_0(u_j^*u_i))$ must be a positive 
$n\times n$ matrix whenever  
$u_1,\dots,u_n\in G$ have height $\leq N$.  
\end{rem}

{\it proof of Lemma \ref{keyLemma}.} We identify the unit $e$ 
of $A$ with its image $\delta_{(0)}\cdot e\in G$.  
Fix $i$, $1\leq i\leq n$,  and write 
$u_ie=\delta_{(n_1,\dots,n_k)}\cdot a_1\otimes\cdots\otimes a_k$.  
Note that $n_k$ must be $0$ because $u_i$ has been multiplied 
on the right by $e$.  

If $n_1>0$ we choose $\ell$, $1<\ell < k$ such that 
$n_1,n_2,\dots,n_\ell$ are positive and $n_{\ell+1}=0$.  
Setting 
$v_i=\delta_{(n_1-1,\dots,n_\ell-1)}\cdot a_1\otimes\cdots\otimes a_\ell$
 and $w_i=\delta_{(n_{\ell+1},\dots,n_k)}
\cdot a_{\ell+1}\otimes\cdots\otimes a_k$, we 
obtain a factorization $u_i=\alpha(v_i)w_i$, 
and we define $b_i$ and $c_i$ by $b_i=E_0(w_i)$, $c_i=0$.  
If $n_1=0$ then $u_ie$ cannot be 
factored in this way; still, we set $v_i=e$, 
and $b_i=c_i=E_0(u_i)$.  This defines 
$b_i$, $c_i$ and $v_i$.  

One now verifies (\ref{keyEq}) in cases:
where both $u_ie$ and $u_je$ factor into a product 
of the form $\alpha(v)w$, when one of 
them so factors and the other does not, and when 
neither does.   For example, if $u_ie=\alpha(v_i)w_i$ 
and $u_je=\alpha(v_j)w_j$ both factor, then we can make 
use of the formulas $E_0(f)=eE_0(f)e=E_0(efe)$ 
for $f\in \ell^1(\Sigma)$, 
$E_0(fg)=E_0(f)E_0(g)$ for $f,g\in [A\ell^1(\Sigma)A]$,
and $E_0\circ\alpha=\phi\circ E_0$, 
to write 
\begin{align*}
E_0(u_j^*u_i)&=E_0(eu_j^*u_ie)=E_0((u_je)^*u_ie)=
E_0(w_j^*\alpha(v_j^*v_i)w_i)\\
&=b_j^*E_0(\alpha(v_j^*v_i))b_i
=b_j^*\phi(E_0(v_j^*v_i))b_i.
\end{align*}

If $u_ie=\alpha(v_i)w_i$ so factors and $u_je=eu_je$ 
does not, then we write
\begin{align*}
E_0(u_j^*u_i)&=E_0((u_je)^*u_ie)=E_0(u_je)^*E_0(\alpha(v_i)w_i)
=b_j^*E_0(\alpha(v_i))b_i\\
&=b_j^*\phi(E_0(v_i))b_i=
b_j^*(\phi(E_0(v_j^*v_i))b_i, 
\end{align*}
noting that in 
this case $v_j^*=e$.  
A similar string of identities settles the 
case $u_ie=eu_ie$, $u_je=\alpha(v_j)w_j$.  

Note that in each of the preceding three cases, 
the terms $c_j^*(e-\phi(e))c_i$ were all $0$.  

In the remaining case where $u_ie=eu_ie$ and 
$u_je=eu_je$, we can write 
$E_0(u_j^*u_i)=E_0(eu_j^*u_ie)=E_0((eu_je)^*eu_ie)=
E_0(eu_je)^*E_0(eu_ie)=b_j^*b_i$.  
Formula (\ref{keyEq}) persists for this case too, 
since $v_i=v_j=e$ and we can write
$$
b_j^*b_i=b_j^*\phi(e)b_i+b_j^*(e-\phi(e))b_i =
b_j^*\phi(E_0(v_j^*v_i))b_i+c_j^*(e-\phi(e))c_i.  \qed
$$

\section{the hierarchy of dilations}\label{S:hierarchy}
Let $(A,\phi)$ be a pair consisting of 
an arbitrary \cstar\ $A$ and a completely positive linear map 
$\phi: A\to A$ 
satisfying $\|\phi\|\leq 1$.  
\begin{defn}\label{dilationDef}
A dilation of $(A,\phi)$ is an $A$-dynamical system 
$(i,B,\alpha)$ with the property that there is an 
$\alpha$-expectation $E: B\to A$ satisfying 
$$
E(\alpha(a))=\phi(a), \qquad a\in A. 
$$
\end{defn}
Notice that the $\alpha$-expectation 
$E: B\to A$ associated with a dilation of 
$(A,\phi)$ is uniquely determined, by
Theorem \ref{universalDilationThm}.  
The class of all dilations of $(A,\phi)$ is contained 
in the class of all $A$-dynamical systems, and it is significant 
that it is also a subcategory.  More explicitly, if 
$(i_1, B_1, \alpha_1)$ and $(i_1, B_2, \alpha_2)$ are two 
dilations of $(A,\phi)$, and if $\theta: B_1\to B_2$ 
is a homomorphism of $A$-dynamical systems, then the 
respective $\alpha$-expectations 
$E_1$, $E_2$ must also transform consistently 
\begin{equation}\label{equivarianceEq}
E_2\circ\theta=E_1.   
\end{equation}
This follows from
Theorem \ref{uniquenessThm}, 
since both $E_1$ and $E_2\circ\theta$ are 
$\alpha_1$-expectations that 
project $\alpha_1(a)$ to $\phi(a)$, $a\in A$.  

Theorem \ref{universalDilationThm} 
implies that every pair $(A,\phi)$ can be 
dilated to the universal $A$-dynamical system 
$(i,{\mathcal P}A,\sigma)$.  Let $E_\phi: {\mathcal P}A\to A$ 
be the $\sigma$-expectation 
satisfying $E_\phi(\sigma(a))=\phi(a)$, 
$a\in A$.  The preceding remarks imply that for every 
other dilation $(i,B,\alpha)$, there is a unique 
surjective $*$-homomorphism $\theta: {\mathcal P}A\to B$ 
such that $\theta\circ\sigma=\alpha\circ\theta$, 
$E\circ\theta=E_\phi$, and which fixes 
$A$ elementwise.  Thus, {\it $(i,{\mathcal P}A,\sigma)$ 
is a universal dilation of $(A,\phi)$}.  The universal 
dilation is obviously too large, since its structure 
bears no relation to $\phi$.  Thus it is significant that 
there is a smallest $(A,\phi)$ dilation, whose structure 
is more closely tied to $\phi$.  We now discuss 
the basic properties of this minimal dilation; 
we examine its structure in 
section \ref{S:corners}.  

In general, every completely positive map 
of $C^*$-algebras $E: B_1\to B_2$ gives rise to a 
norm-closed two-sided ideal $\ker E$ in $B_1$ as follows
\begin{equation*}
\ker E=\{x\in B_1: E(bxc)=0, \quad {\text {for all }} b,c\in B_1 \}.  
\end{equation*}
In more concrete terms, 
if $B_2\subseteq{\mathcal B}(H)$ acts concretely 
on some Hilbert space and $E(x)=V^*\pi(x)V$ is a Stinespring 
decomposition of $E$, where $\pi$ is a representation of 
$B_1$ on some other Hilbert space $K$ and $V: H\to K$ is a bounded 
operator such that $\pi(B_1)VH$ has $K$ as its 
closed linear span, then one 
can verify that 
\begin{equation}\label{kernelEq}
\ker E=\{x\in B_1: \pi(x)=0\}.
\end{equation}
Notice too that $\ker E=\{0\}$ iff $E$ is {\it faithful on ideals}
in the sense that for every two-sided ideal $J\subseteq B$, 
one has $E(J)=\{0\}\implies J=\{0\}$.  
\begin{prop}\label{kernelProp}
Let $(i,B,\alpha)$ be an $A$-dynamical system and 
let $E: B\to A$ be an $\alpha$-expectation.  Then 
$\ker E$ is an $\alpha$-invariant ideal with the property 
$A\cap\ker E=\{0\}$.  

If $(i_1,B_1,\alpha_1)$ and $(i_2,B_2,\alpha_2)$ 
are two dilations of $(A,\phi)$ and $\theta: B_1\to B_2$ 
is a homomorphism of $A$-dynamical systems, then 
\begin{equation}\label{expectationEq}
\ker E_1=\{x\in B_1: \theta(x)\in\ker E_2\}.
\end{equation}
\end{prop}

\begin{Prf}
That $A\cap \ker E=\{0\}$ is clear 
from the fact that if $a\in A\cap \ker E$ then 
$AaA=E(AaA)=\{0\}$, hence $a=0$.  Relation (\ref{expectationEq})
is also straightforward, since $\theta(B_1)=B_2$ and 
$E_2\circ\theta=E_1$.  Indeed, for each $x\in B_1$, we have 
$$
E_2(B_2\theta(x)B_2)=E_2(\theta(B_1)\theta(x)\theta(B_1))=
E_2(\theta(B_1xB_1))=E_1(B_1xB_1),
$$
from which (\ref{expectationEq}) follows.

To see that $\alpha(\ker E)\subseteq \ker E$, 
choose $k\in \ker E$.  Since $B$ is spanned by all finite 
products of elements $\alpha^n(a)$, $a\in A$, 
$n=0,1,\dots$, it suffices to show 
that $E(y\alpha(k)x)=0$ for all $y\in B$ and all $x$ 
of the form
$x=\alpha^{m_1}(a_1)\cdots\alpha^{m_k}(a_k)$.
Being a completely positive contraction, 
$E$ satisfies the Schwarz inequality
$$
E(y\alpha(k)x)^*E(y\alpha(k)x)\leq 
E(x^*\alpha(k^*)y^*y\alpha(k)x)\leq 
\|y\|^2E(x^*\alpha(k^*k)x);
$$
hence it suffices to show that $E(x^*\alpha(k^*k)x)=0$.  
To prove the latter, one can argue cases as follows.  
Assuming that $m_i>0$ for all $i$, then 
$x=\alpha(x_0)$ for some
$x_0\in B$, and using $E\circ\alpha=\phi\circ E$ 
one has 
$$
E(x^*\alpha(k^*k)x)=E(\alpha(x_0k^*kx_0))=
\phi(E(x_0^*k^*kx_0))=0.
$$
For the remaining case where some $m_i=0$,
notice 
that $x$ must have one of the forms $x=a\in A$ 
(when all $m_i$ are $0$), or $x=ax_0$ with $x_0\in B$ 
(when $m_1=0$ and some other $m_j$ is positive), 
or $\alpha(x_1)ax_2$ 
with $a\in A$, $x_1,x_2\in B$ 
(when $m_1>0, \dots, m_{r-1}>0$ and $m_r=0$), 
and in each case $E(x^*\alpha(k^*k)x)=0$.   
If $x=\alpha(x_1)ax_2$, for example, then 
we make use of Formula (\ref{multipEq}) to write
$$
E(x^*\alpha(k^*k)x)=E(x_2^*a^*\alpha(x_1^*k^*kx_1)ax_2)=
E(x_2)^*a^*E(\alpha(x_1^*k^*kx_1)aE(x_2).
$$
The term on the right vanishes 
since $E(\alpha(x_1^*k^*kx_1)=\phi(E(x_1^*k^*kx_1))=0$.  
The other cases are dealt with similarly, 
and $\alpha(\ker E)\subseteq \ker E$ follows.   
\qed
\end{Prf}
\vskip0.1in
We deduce the existence of minimal dilations and 
their basic characterization as follows.  
Fix a pair $(A,\phi)$, where 
$\phi: A\to A$ is a completely positive contraction, 
let $\sigma$ be the shift on ${\mathcal P}A$, 
and let $E_\phi: {\mathcal P}A\to A$ 
be the unique $\sigma$-expectation satisfying 
$E_\phi(\sigma(a))=\phi(a)$, 
$a\in A$.  Proposition \ref{kernelProp} implies that  
$\sigma$ leaves $\ker E_\phi$ 
invariant, thus it can be promoted to an 
endomorphism $\dot\sigma$ of the quotient $C^*$-algebra 
${\mathcal P}A/\ker E_\phi$.  Moreover, since 
$A\cap\ker E_\phi=\{0\}$, the inclusion of $A$ in 
${\mathcal P}A$ promotes to an inclusion of 
$A$ in ${\mathcal P}A/\ker E_\phi$.  Thus we obtain an 
$A$-dynamical system $(i,{\mathcal P}A/\ker E_\phi,\dot\sigma)$
having a natural 
$\dot\sigma$-expectation $\dot E$ defined by 
$\dot E(x+\ker E_\phi)=E_\phi(x)+\ker E_\phi$, which satisfies 
$\dot E(\dot\sigma(a))=\phi(a)$, $a\in A$.  
It is called 
the {\it minimal} dilation of $(A,\phi)$ in light    
of the following:

\begin{cor}\label{minimalCor}
The dilation $(i,{\mathcal P}A/\ker E_\phi,\dot\sigma)$ 
of $(A,\phi)$ has the following properties.  
\begin{enumerate}
\item
$(i,{\mathcal P}A/\ker E_\phi,\dot\sigma)$ is subordinate to all 
other dilations of $(A,\phi)$.  
\item
The $\dot\sigma$-expectation $\dot E$ of 
$(i,{\mathcal P}A/\ker E_\phi,\dot\sigma)$ 
satisfies $\ker\dot E=\{0\}$.  
\item
Every dilation $(i,B,\alpha)$ of $(A,\phi)$ whose 
$\alpha$-expectation $E$ satisfies $\ker E=\{0\}$ is isomorphic 
to $(i,{\mathcal P}A/\ker E_\phi,\dot\sigma)$.  
\end{enumerate}
\end{cor}

\begin{Prf}
(2) follows by construction 
of $(i,{\mathcal P}A/\ker E_\phi,\dot\sigma)$, 
since the kernel ideal of 
its $\dot\sigma$-expectation has been reduced to $\{0\}$.  

To prove (1), let $(i,B,\alpha)$ be an arbitrary 
dilation of $(A,\phi)$.  By the universal property 
of $(i,{\mathcal P}A,\sigma)$ there is a surjective 
$*$-homomorphism $\theta: {\mathcal P}A\to B$ satisfying 
$\theta\circ\sigma=\alpha\circ\theta$; and by 
(\ref{equivarianceEq}) one has $E\circ\theta=E_\phi$.  
Formula (\ref{expectationEq}) implies that 
$\ker E_\phi$ contains $\ker\theta$, hence we 
can define a morphism of $C^*$-algebras 
$\omega:B \to {\mathcal P}A/\ker E_\phi$ 
by way of $\omega(\theta(x))=x+\ker E_\phi$, 
for all $x\in {\mathcal P}A$.  Obviously, 
$\omega$ is a homomorphism of $A$-dynamical 
systems, and we conclude that 
$(i,B,\alpha)\geq (i,{\mathcal P}A/\ker E_\phi,\dot\sigma)$.  

For (3), notice that if $(i,B,\alpha)$ is a dilation 
of $(A,\phi)$ whose  $\alpha$-expectation $E: B\to A$ 
satisfies $\ker E=\{0\}$ 
and $\theta: {\mathcal P}A\to B$ is the homomorphism 
of the previous paragraph, then Formula (\ref{expectationEq}) 
implies that $\ker E_\phi=\ker\theta$.  Thus 
$\omega:B \to {\mathcal P}A/\ker E_\phi$ 
has trivial kernel, hence it must implement 
an isomorphism of $A$-dynamical systems 
$(i,B,\alpha)\cong(i,{\mathcal P}A/\ker E_\phi,\dot\sigma)$.\qed
\end{Prf}

\section{structure of minimal dilations}\label{S:corners}

Corollary \ref{minimalCor}
implies that minimal dilations of $(A,\phi)$ exist for 
every contractive completely positive map $\phi: A\to A$, 
and that they are characterized by the fact that their 
$\alpha$-expectations are faithful on ideals.  The 
latter imposes strong requirements  
on the structure of minimal dilations, and we conclude 
by elaborating on these structural issues.

\begin{defn}\label{standardDef}
A {\it standard} dilation of $(A,\phi)$ is a dilation $(i,B,\alpha)$ 
such that $A=pBp$ is an essential corner of $B$ 
whose projection $p\in M(B)$ satisfies 
$p\alpha(x)p=\phi(pxp)$, $x\in B$. 
\end{defn}
 In such cases, 
$E(x)=pxp$ is the $\alpha$-expectation 
of $B$ on $A$.  
Standard dilations are most transparent 
in the special case where 
$A$ has a unit $e$ and $\phi(e)=e$.  
To illustrate that, let $B$ be a \cstar\ containing $A$ and 
let $\alpha$ be an endomorphism of $B$ with the following 
property: 
\begin{equation}\label{stdEq}
\phi(a)=e\alpha(a)e,\qquad a\in A, \label{oldEq6}
\end{equation}
$e$ denoting the unit of $A$.  We may also assume 
that $B$ is generated 
by $A\cup \alpha(A)\cup\alpha^2(A)\cup\cdots$, 
so that $(i,B,\alpha)$ becomes an 
$A$-dynamical system.  

\begin{prop}\label{standardProp}
The projection $e\in B$ satisfies $\alpha(e)\geq e$, 
$A=eBe$ is a hereditary subalgebra 
of $B$, and the map $E(x)=exe$ defines an $\alpha$-expectation
from $B$ to $A$.  If, in addition, $A$ is an essential subalgebra 
of $B$, then $(i,B,\alpha)$ 
is a standard dilation of $(A,\phi)$.  

\end{prop}

{\it Sketch of proof.}
Formula (\ref{stdEq}) implies that $e\alpha(e)e=\phi(e)=e$, 
hence $\alpha(e)\geq e$.  It follows immediately that  
$e\alpha(exe)e=e\alpha(x)e$ for 
$x\in B$.  

At this point, a simple induction establishes 
$e\alpha^n(a)e=\phi^n(a)$, $a\in A$, $n=0,1,2,\dots$.  
An argument similar to the proof 
of Theorem \ref{uniquenessThm} allows one to evaluate 
more general expectation values as in (\ref{nPointEq})
$$
e\alpha^{n_1}(a_1)\alpha^{n_2}(a_2)\cdots\alpha^{n_k}(a_k)e 
=[n_1,\dots,n_k; a_1,\dots,a_k],
$$
which implies $eBe\subseteq A$.  Hence
$A=eBe$ is a hereditary subalgebra of $B$.  

With these formulas in hand one finds 
that the conditional expectation 
$E(x)=exe$ satisfies axioms E1 and E2
of Definition \ref{dilationDef}.  Hence $(i,B,\alpha)$ 
is a standard dilation of $(A,\phi)$ whenever 
$A$ is an essential subalgebra of $B$.  \qed
\vskip0.1in
We remark that the converse is also true: 
given $(A,\phi)$ for which $A$ has a unit $e$ 
and $\phi(e)=e$, then every standard dilation 
has the properties 
of Proposition (\ref{standardProp}).  
The description of standard dilations in 
general, where $A$ is 
unital and $\|\phi\|\leq 1$ or is perhaps nonunital, 
becomes more subtle.  

The universal dilation $(i,{\mathcal P}A,\sigma)$
of $(A,\phi)$ is {\it not} a standard  dilation.  
For example, when $A$ has a unit $e$ one can make use of 
the universal property of 
${\mathcal P}A=A*A*\cdots$ to 
exhibit representations $\pi:{\mathcal P}A\to {\mathcal B}(H)$ 
such that $\pi(e)$ and 
$\pi(\sigma(e))$ are nontrivial orthogonal projections.   
Hence $\sigma(e)\ngeq e$.  
Moreover, $A$ is not a hereditary subalgebra of ${\mathcal P}A$, 
and the conditional expectation 
of Theorem \ref{universalDilationThm} is  
never of the form $x\mapsto exe$.  

On the other hand, we now show 
that minimal dilations of $(A,\phi)$ must be standard.  
This is based on the following characterization of 
essential corners in terms of conditional expectations.

\begin{prop}\label{essCornerProp}
For every inclusion of $C^*$-algebras $A\subseteq B$, 
the following are equivalent.  
\begin{enumerate}
\item[(i)]
$A$ is an essential corner $pBp$ of $B$.
\item[(ii)]
There is a conditional expectation $E: B\to A$ 
whose restriction to $[ABA]$ is multiplicative, and
which satisfies $\ker E=\{0\}$.
\end{enumerate}
Moreover, the conditional 
expectation $E:B\to A$ of (ii) is 
the compression map $E(x)=pxp$, and 
it is unique.  The projection $p\in M(B)$ 
satisfies 
$$
\lim_{n\to\infty}\|xe_n-xp\|=0,\qquad x\in B,
$$
where $e_n$ is any approximate unit for $A$, and it 
defines the closed left ideal generated by $A$ as 
follows: $[BA]=Bp$.  
\end{prop}

\begin{Prf} The implication (i) $\implies$ (ii) is 
straightforward, since the compression 
map $E(x)=pxp$ obviously defines a 
conditional expectation of $B$ on $A=pBp$ that 
is multiplicative on $[ABA]$.  If $x\in B$ satisfies 
$E(BxB)=\{0\}$ then $pBx^*xBp=\{0\}$, hence 
$xBA=xBpBp=\{0\}$, and therefore $x=0$ 
because $[BAB]$ is assumed to be an essential ideal in (i).   

(ii) $\implies$ (i).  Given a conditional expectation 
$E: B\to A$ satisfying (ii), we may assume that 
$A\subseteq{\mathcal B}(H)$ acts nondegenerately on 
some Hilbert space (e.g., represent $B$ faithfully 
on some Hilbert space and take $H$ to be the closed 
linear span of the ranges of all operators in $A$).  
Thus $E: B\to {\mathcal B}(H)$ becomes an operator-valued 
completely positive map of norm $1$, having 
a Stinespring decomposition $E(x)=V^*\pi(x) V$,
with $\pi$ a representation of $B$ on a Hilbert space 
$K$, and $V: H\to K$ a contraction 
with $[\pi(B)VH]=K$. 

Let $P$ be the projection on $[\pi(A)K]$.  
We claim that $V$ is an 
isometry with $VV^*=P$.  To prove that, choose $a\in A$, 
$b\in B$, and let $e_n$ be an approximate unit for $A$.  
Since $E$ is multiplicative on $[ABA]$ we can write 
\begin{align*}
V^*\pi(b^*a^*)(VV^*-{\bf 1})\pi(ab)V^* =
E(b^*a)E(ab) -E(b^*a^*ab)&=\\
\lim_ne_n(E(b^*a^*)E(ab)-E(b^*a^*ab))e_n &=\\
\lim_n(E(e_nb^*a^*)E(abe_n)-E(e_nb^*a^*abe_n)) &= 0.  
\end{align*}
It follows that $VV^*-{\bf 1}$ vanishes on the closed 
linear span of $\pi(A)\pi(B)VH$, namely $[\pi(A)K]$; 
hence $VV^*\geq P$.  
On the other hand,
for $a\in A$ we have 
$Va=VE(a)=VV^*\pi(a)V=\pi(a)V$. Thus 
$VH\subseteq [VAH]=[\pi(A)VH]\subseteq PK$; hence  $VV^*\leq P$.  
That $V$ is an isometry follows from the fact that
 for $a\in A$, 
$V^*Va=V^*VE(a)=V^*VV^*\pi(a)V=V^*\pi(a)V=E(a)=a$, 
and by nondegeneracy $H$ is the closed 
linear span of $\{a\xi: a\in A,\quad \xi\in H\}$.

We claim that $P=VV^*$ belongs to the multiplier algebra 
of $\pi(B)$.  For that, choose an approximate 
unit $e_n$ for $A$.  Since both 
$\pi(e_n)$ and $VV^*$ are self-adjoint, 
it suffices to show that 
for every $b\in B$, $\pi(b)\pi(e_n)\to \pi(b)VV^*$ in 
norm as $n\to\infty$.  Using $VV^*\pi(e_n)=\pi(e_n)VV^*=\pi(e_n)$, 
we can write 
\begin{align*}
&\|\pi(b)(\pi(e_n)-VV^*)\|^2=
\|(\pi(e_n)-VV^*)\pi(b^*b)(\pi(e_n)-VV^*)\|=\\
&\|VV^*(\pi(e_nb^*be_n)-\pi(e_nb^*b)-\pi(b^*be_n)+\pi(b^*b))VV^*\|=\\
&\|V(E(e_nb^*be_n)-E(e_nb^*b)-E(b^*be_n)+E(b^*b) )V^* \|\leq \\
&\|e_nE(b^*b)e_n-e_nE(b^*b)-E(b^*b)e_n+E(b^*b)\|,
\end{align*}
and the last term tends to $0$ as $n\to\infty$ 
because $e_n$ is an approximate unit for $A$  
and $E(b^*b)\in A$.  It 
follows that $\pi(B)VV^*=\pi(A)P$ is the closed left ideal 
in $\pi(B)$ generated by $\pi(A)$.  

We claim next that $\pi(A)=P\pi(B)P$ is a corner of $\pi(B)$.  
Indeed, 
$$
VV^*\pi(B)VV^*=VE(B)V^*=\pi(E(B))VV^*=\pi(A)VV^*=\pi(A).
$$  
It is essential because for any operator
$T\in{\mathcal B}(K)$ for which $T\pi(BA)=\{0\}$ we must have 
$T\pi(BA)VV^*=\{0\}$.  But since $K$ is spanned 
by vectors of the form $\pi(b)Va\xi=\pi(b)\pi(a)V\xi$ 
for $a\in A$, $b\in B$, the only possibility is $T=0$.  

Finally, $\pi$ must be a faithful representation because 
$\pi(x)=0$ implies
$$
E(BxB)=V^*\pi(B)\pi(x)\pi(B)V=\{0\},
$$ 
and the latter implies $x=0$ by hypothesis (ii).  The 
preceding assertions can now be pulled back through 
the isomorphism $\pi: B\to \pi(B)$ to give (i).  
\qed
\end{Prf}
\vskip0.1in
Combining Corollary \ref{minimalCor} with 
Proposition \ref{essCornerProp}, 
we obtain:  

\begin{thm}\label{cornerThm}
For every $(A,\phi)$ as above, the minimal 
dilation of $(A,\phi)$ is a standard dilation 
satisfying the assertions of Proposition \ref{essCornerProp}.  
All standard dilations of $(A,\phi)$ are equivalent to 
the minimal one.  
\end{thm}

We remark that Theorem \ref{cornerThm}, 
together with a theorem of Larry Brown \cite{lgbStable}, implies 
that the \cstar\ $B$ of the {\it minimal} dilation $(i,B,\alpha)$ 
of $(A,\phi)$ can be embedded in the multiplier 
algebra of $A\otimes\mathcal K$.

\subsection*{Concluding Remarks}
It is appropriate to review some highlights of the literature 
on noncommutative dilation theory, since it bears 
some relationship to the contents of \S\S\ref{S:hierarchy}--\ref{S:corners}.  
Several approaches to dilation theory for semigroups 
of completely positive maps have been proposed since 
the mid 1970s, including work of Evans and Lewis \cite{EvLewis}, 
Accardi et al \cite{AFL}, K\"ummerer \cite{kumMarkovDil}, 
Sauvageot \cite{sauv86}, 
and many others.  Our attention was drawn to these 
developments by work of Bhat and Parthasarathy 
\cite{BhatPar}, in which the first dilation 
theory for CP semigroups acting on ${\mathcal B}(H)$ 
emerged that was effective for our work on $E_0$-semigroups
\cite{arvPureAbsorbing}, \cite{arvInteractions}.  
SeLegue \cite{selThesis} 
showed how to apply multi-operator dilation theory 
to obtain the Bhat--Parthasarathy results, and 
he calculated the expectation values of the 
$n$-point functions of such dilations.  Recently, 
Bhat and Skeide \cite{BhSkeide} have initiated an 
approach to the subject 
that is based on Hilbert modules over \cstar s 
and von Neumann algebras.  

We intend to take up applications to semigroups of 
completely positive maps elsewhere.

\vfill
\eject

\bibliographystyle{alpha}


\begin{thebibliography}{K{\"u}m85}

\bibitem[AL82]{AFL}
A.~Accardi, L.~Frigerio and J.~T. Lewis.
\newblock Quantum stochastic processes.
\newblock {\em Publ. RIMS, Kyoto Univ.}, 18:97--133, 1982.

\bibitem[Arv97]{arvPureAbsorbing}
William Arveson.
\newblock Pure {$E_0$}-semigroups and absorbing states.
\newblock {\em Comm. Math. Phys.}, 187:19--43, (1997).

\bibitem[Arv00]{arvInteractions}
William Arveson.
\newblock Interactions in noncommutative dynamics.
\newblock {\em Comm. Math. Phys.}, 211:63--83, (2000).

\bibitem[BP94]{BhatPar}
B.~V.~R. Bhat and K.~R. Parthasarathy.
\newblock Kolmogorov's existence theorem for markov processes in
  {$C^*$}-algebras.
\newblock {\em Proc. Indian Acad. Sci. (Math. Sci.)}, 104:253--262, 1994.

\bibitem[Bro77]{lgbStable}
L.~G. Brown.
\newblock Stable isomorphism of hereditary subalgebras of {$C^*$}-algebras.
\newblock {\em Pac. J. Math.}, 71(2):335--348, 1977.

\bibitem[BS00]{BhSkeide}
B.~V.~R. Bhat and Michael Skeide.
\newblock Tensor product systems of hilbert modules and dilations of completely
  positive semigroups.
\newblock In {\em Infinite dimensional analysis, quantum probability and
  related topics}, volume~3, pages 519--575, 2000.

\bibitem[EL77]{EvLewis}
D.~E. Evans and J.~T. Lewis.
\newblock Dilations of irreversible evolutions in algebraic quantum theory.
\newblock {\em Comm. Dublin Inst. Adv. Studies Series A24}, page what pages?,
  1977.

\bibitem[K{\"u}m85]{kumMarkovDil}
B.~K{\"u}mmerer.
\newblock Markov dilations on $w^*$-algebras.
\newblock {\em J. Funct. Anal.}, 63:139--177, 1985.

\bibitem[Ped79]{gPedAutbook}
G.~K. Pedersen.
\newblock {\em {$C^*$}-algebras and their automorphism groups}.
\newblock Academic Press, London, 1979.

\bibitem[Sau86]{sauv86}
J.-L. Sauvageot.
\newblock Markov quantum semigroups admit covariant {$C^*$}-dilations.
\newblock {\em Comm. Math. Phys.}, 106:91--103, 1986.

\bibitem[SeL97]{selThesis}
Dylan SeLegue.
\newblock {\em Minimal Dilations of CP Maps and a {$C^*$}-Extension of the
  Szeg\"o Limit}.
\newblock PhD thesis, University of California, Berkeley, June 1997.

\end{thebibliography}
\newcommand{\noopsort}[1]{} \newcommand{\printfirst}[2]{#1}
  \newcommand{\singleletter}[1]{#1} \newcommand{\switchargs}[2]{#2#1}

\end{document}